\documentclass{article}

\usepackage{arxiv}

\usepackage[utf8]{inputenc} 
\usepackage[T1]{fontenc}    
\usepackage{hyperref}       
\usepackage{url}            
\usepackage{booktabs}       
\usepackage{amsfonts}       
\usepackage{nicefrac}       
\usepackage{microtype}      
\usepackage{lipsum}
\usepackage{amsmath,amsthm}
\usepackage{amssymb}
\newtheorem{theorem}{Theorem}

\newtheorem{lemma}[theorem]{Lemma}

\newtheorem{remark}[theorem]{Remark}

\newcommand{\beq} {\begin{eqnarray*}}
\newcommand{\eeq} {\end{eqnarray*}}

\title{Exact rate of convergence of the mean Wasserstein distance between the empirical and true Gaussian distribution}

\newcommand{\F}{\mathbb{F}}
\newcommand{\G}{\mathbb{G}}
\newcommand{\E}{\mathbb{E}}

\newcommand{\Var}{\hbox{{\rm Var}}}

\newcommand{\Cov}{\hbox{{\rm Cov}}}

\author{Philippe Berthet\thanks{Institut de Math\'{e}matiques de Toulouse UMR 5219 ; Universit\'{e} Paul Sabatier, France. \it{philippe.berthet@math.univ-toulouse.fr}}
   \And
 Jean Claude Fort\thanks{MAP5 UMR 8145; Universit\'e Paris-Descartes,
    France. \it{jean-claude.fort@parisdescartes.fr}}}

\begin{document}
\maketitle
\begin{abstract}
We study the Wasserstein distance $W_2$ for Gaussian samples. We establish the exact rate of convergence $\sqrt{\log\log n/n}$ of the expected value of the $W_2$ distance between the empirical and true $c.d.f.$'s for the normal distribution. We also show that the rate of weak convergence is unexpectedly $1/\sqrt{n}$ in the case of two correlated Gaussian samples.
\end{abstract}

{\bf AMS} Subject Classification: {62G30 ; 62G20 ; 60F05 ; 60F17}

\section{Introduction}\label{sec:intro}

In this article we investigate in details the asymptotic behaviour of the quadratic Wasserstein distance between the empirical cumulative distribution function ($c.d.f.$) of a sample $X_1,\dots ,X_n$ of independent standard Gaussian random variables denoted by $\mathbb F_n$ and the standard normal $c.d.f.$ denoted by $\Phi$. Thus we consider the random variable
\[W_2^2(\mathbb F_n,\Phi)=\int_0^1\vert\mathbb F_n^{-1}(u)-\Phi^{-1}(u)\vert^2du.
\]
More precisely we are interested in the exact rate of convergence of $\mathbb E \left(W_2^2(\mathbb F_n,\Phi)\right)$. Define $h(u)=\Phi'\circ\Phi^{-1}(u)$ for $u\in(0,1)$. First note that Corollary 19 in \cite{BF19} does not apply in this specific case where $b=2$, and indeed we almost surely have $\lim_{n\rightarrow+\infty} n W_2^2(\mathbb F_n,\Phi) = + \infty$. Secondly, to our knowledge the most precise result about the behaviour of $W_2(\mathbb F_n,\Phi)$ is given by Theorem 4.6 (ii) in \cite{DelBarrio05} which implies, as $n\rightarrow +\infty$, the convergence in distribution
\begin{equation}\label{Del}
n W_2^2(\mathbb F_n,\Phi)-\int_{1/n}^{1-1/n}\frac{u(1-u)}{h^{2}(u)}du\to
{\displaystyle\int\nolimits_{0}^{1}}
\frac{\mathbb{B}^{2}(u)-\mathbb{E}\left(  \mathbb{B}^{2}(u)\right)  }
{h^{2}(u)}du,
\end{equation}
where $\mathbb B$ is a standard Brownian bridge. This is not enough to control $n\mathbb E (W_2^2(\mathbb F_n,\Phi))$ since the deterministic centering integral is diverging. In \cite{Ledoux17} specific bounds on $n\mathbb E ( W_p^p(\mathbb F_n,F))$ are given for $\log$-concave distribution $F$. In the standard Gaussian case Corollary 6.14 of \cite{Ledoux17} reads
\begin{equation}\label{BL}
c\frac{\log \log n}{n} \leqslant \mathbb E \left(W_2^2(\mathbb F_n,\Phi)\right) \leqslant C\frac{\log \log n}{n}
\end{equation}
where $0<c<C<+\infty $. The main achievement below is to compute the exact asymptotic constant in (\ref{BL}). As far as we know this is the first result of this kind.

\noindent In the spirit of \cite{BF19} we moreover extend the investigations in the one sample case to the two correlated samples case. More precisely, we study the random quantity $ W_2^2(\mathbb \F_n,\G_n)$ where $\F_n,\G_n$ are the marginal empirical $c.d.f.$ obtained from a $n$-sample $(X_i,Y_i)_{ 1\leqslant i\leqslant n}$ of standard Gaussian couples with correlation $\rho$. If the Gaussian marginals $\Phi_X$ and $\Phi_Y$ were not identical the general Theorem 14 in \cite{BFK17} would imply the convergence in distribution
\begin{equation}\label{BFK}
\sqrt{n} \bigl( W_{2}^2(\mathbb{F}_{n},\mathbb{G}_{n})-W^2_{2}(\Phi_X,\Phi_Y)
\bigr) \rightarrow \mathcal{N}\bigl(0,\sigma ^{2}(\Sigma)
\bigr)
\end{equation}
where $\Sigma$ is the covariance matrix of $(X_1,Y_1)$ and $\sigma ^{2}(\Sigma)$ has a closed form expression that explicitly depends on $\Sigma$. In particular, Corollary 18 of \cite{BFK17} shows that for two independent samples from two distinct Gaussian distributions $\mathcal{N}(\nu ,\zeta ^{2})$ and $\mathcal{N}(\mu ,\xi ^{2})$ it holds $\sigma ^{2}(\Sigma)=4(\zeta ^{2}+ \xi ^{2})(
\nu -\mu )^{2}+2(\zeta ^{2}+ \xi ^{2})(\zeta - \xi )^{2} $.

\noindent Surprisingly, the second result below establishes that whenever the marginals are the same, $\Phi_X=\Phi_Y=\Phi$, and the samples are not independent, that is $\rho\ne0$, the rate of weak convergence of $\ W_2^2(\mathbb \F_n,\G_n)$ is  $1/n$ and the limiting distribution is a slight variation of the one given at Theorem 11 in \cite{BF19}, even if the sufficient condition of the latter result is not satisfied.

\section{The results}

First we provide the limiting constant in \eqref{BL}.
\begin{theorem}\label{CVE}
Let $\mathbb F_n$ be the empirical $c.d.f.$ of an $i.i.d.$ standard normal sample of size $n$ and $\Phi$ the  $c.d.f.$ of the standard normal distribution. Then it holds
\begin{align*}
\lim_{n\rightarrow+\infty}\frac{n}{\log\log n}\mathbb{E}\left(W_{2}%
^{2}(\mathbb F_n,\Phi)\right)   &  =1,\\
\lim_{n\rightarrow+\infty}\sqrt{\frac{n}{\log\log n}}\mathbb{E}\left( 
W_{2}(\mathbb F_n,\Phi)\right)   &  =1.
\end{align*}
\end{theorem}
\begin{remark}\label{Rem0}
This result is consistent with (\ref{Del}) and the fact that, by \cite{Bickel78}, we have
$${\displaystyle\int\nolimits_{1/n}^{1-1/n}}
\frac{u(1-u)}{h^{2}(u)}du=\log \log n + \log 2 + \gamma _0 + o(1)$$
which implies that $\frac{n}{\log\log n}W_{2}
^{2}(\mathbb F_n,\Phi) \to 1$ in probability.
\end{remark}
\begin{remark}\label{Rem1}
In the case of a sample of unstandardized normal random variables with variance $\sigma^2$ the expected $W_2$-distance between the empirical and the true distribution has the same rate as above and limiting constants $\sigma^2$ and $\sigma$, respectively.
\end{remark}
\begin{remark}\label{Remindep}
If $\mathbb G_n$ is a second empirical $c.d.f.$ independent of $\mathbb F_n$ and build from another sample we see that $\mathbb{E}\left(W_{2}%
^{2}(\mathbb F_n,\mathbb G_n)\right)=\mathbb{E}\left(W_{2}%
^{2}(\mathbb F_n,\Phi)\right)+\mathbb{E}\left(W_{2}%
^{2}(\mathbb G_n,\Phi)\right)$ since $\mathbb{E}(
\int_0^1(\mathbb F_n^{-1}(u)-\Phi^{-1}(u))du)=0$. Therefore, in this independent case we have
\begin{align*}
\lim_{n\rightarrow+\infty}\frac{n}{\log\log n}\mathbb{E}\left(W_{2}%
^{2}(\mathbb F_n,\mathbb G_n)\right)   &  =2
\end{align*}
which is in contrast with the forthcoming dependent sample case.
\end{remark}

\noindent Second, in the setting of \cite{BFK17} and \cite{BF19} we also get the rate of weak convergence in the two correlated samples case. 
\begin{theorem}\label{Theo2}
 Let $\mathbb F_n$ and $\mathbb G_n$ denote the marginal empirical $c.d.f.$ of a size $n$ $i.i.d.$ sample of correlated  bivariate standard normal with covariance $\rho$, $0< |\rho |< 1$. Let
\begin{align*}
C_\rho(u,v)&=\mathbb P(X\leqslant \Phi^{-1}(u),Y\leqslant \Phi^{-1}(v)),\quad u,v\in(0,1),\\
\displaystyle {\cal G}(u)&=\frac{\mathbb{B}^{X}(u)}{h(u)}-\frac{\mathbb{B}^{Y}(u)}{h(u)},\quad u\in(0,1),
\end{align*}
where $(\mathbb{B}^{X},\mathbb{B}^{Y})$ are two standard Brownian bridges with cross covariance $$\Cov(\mathbb{B}^{X}(u),\mathbb{B}^{Y}(v))=
C_\rho(u,v)-uv,\  u,v\in\left(
0,1\right).$$
Then we have the convergence in distribution
$$\displaystyle n W_{2}^2(\mathbb F_n,\mathbb G_n)\to ||\mathcal{G}||_2^2=\int_0^1  {\cal G}(u)^2 du$$and the limiting random variable is almost surely finite with finite expectation.
\end{theorem}
\begin{remark}\label{Rem2}
By Theorem \ref{Theo2} it holds $\displaystyle \sqrt{n} W_{2}(\mathbb F_n,\mathbb G_n)\to ||\mathcal{G}||_2$ with a CLT rate and a non degenerate limiting distribution with finite variance. This was not expected since in the case of two independent samples, that is $\rho=0$, it holds
$$
\mathbb{E}(||\mathcal{G}||_2^2)=\int_0^1  \mathbb{E}({\cal G}(u)^2) du=2\int_0^1  \frac{u(1-u)}{h^2(u)}  du=+\infty
$$
which proves by Theorem 1.3 of \cite{CHS93} that $\mathbb{P}(||\mathcal{G}||_2=+\infty)=1$, and is consistent with the similar case where $\mathbb G_n$ is replaced with $\Phi$ as shown by Theorem \ref{CVE}. 
\end{remark}

\begin{remark}\label{Rem3}
Theorem \ref{Theo2} is an extension of Theorem 11 in \cite{BF19} for Gaussian correlated samples that proves that the dependency between two $i.i.d.$ samples expressed through the joint law may influence the rate of convergence of $W_{2}^2(\mathbb F_n,\mathbb G_n)$ if the marginal distributions are the same. In the general CLT formulated at Theorem 14 of \cite{BFK17}, only the limiting finite variance of $\sqrt{n}(W_{2}^2(\mathbb F_n,\mathbb G_n)-W_{2}^2(\Phi_X,\Phi_Y))$ was affected by the joint law if the marginal distributions are different, not the rate $1/\sqrt{n}$ as recalled at \eqref{BFK} above.
\end{remark}
\section{Proofs}\label{sec:proofs}

\subsection{Preliminaries}\label{sec:prelim}

\noindent  Note that the density quantile function $h(u)=\Phi'\circ\Phi^{-1}(u)$ is symmetric on $(0,1)$ about $u=1/2$. Straightforward computations yield, as $x\rightarrow+\infty$, %
\begin{align*}
\psi(x)  &  =-\log(1-\Phi(x))=\frac{x^{2}}{2}+\log x+\frac{1}{2}\log
(2\pi)+O\left(  \frac{1}{x^{2}}\right),\\
\psi^{-1}(x)  &  =\sqrt{2\left(  x-\frac{1}{2}\log x-\frac{1}{2}\log
(2\pi)-\frac{1}{2}\log2+O\left(  \frac{\log x}{x}\right)  \right)  }.
\end{align*}
As a consequence, we have, as $u\rightarrow 1$,
\begin{align}
&\Phi^{-1}(u)    =\psi^{-1}\left(  \log\left(  \frac{1}{1-u}\right)  \right)\nonumber
\\
 &=\sqrt{2\left(  \log\left(  \frac{1}{1-u}\right)  -\frac{1}{2}\log
\log\left(  \frac{1}{1-u}\right)  -\frac{1}{2}\log(4\pi)+O\left(  \frac
{\log\log\left(  1/(1-u)\right)  }{\log\left(  1/(1-u)\right)  }\right)
\right)  },\label{phiinv}
\end{align}
and
\begin{eqnarray}\label{h}
h(u)= \Phi'\circ\Phi^{-1}(u)=\sqrt{2}(1-u)\sqrt{\log\left(
\frac{1}{1-u}\right)  }\left(1+O\left( \frac
{\log\log\left(  1/(1-u)\right)  }{\log\left(  1/(1-u)\right)  }\right)\right).\label{h}
\end{eqnarray}


\noindent Let us extend the results concerning the first and second moments of the extreme order statistics of a Gaussian sample stated at page 376 in \cite{Cramer}.

\begin{lemma}
\label{Cramer} Let $Z_1 \leqslant \dots \leqslant Z_n$ denote the order statistics of $X_1,...,X_n$. Let $1 \leqslant\theta \leqslant 2$ and $C>0$. For any $k\leqslant C(\log n)^\theta$ it holds
\begin{align*}
    \mathbb{E}\left(Z_{n-k} \right)&=\sqrt{2\log n}-\frac{\log\log n+2(s_{k+1}^1-\gamma_0)+\log(4\pi)}{\sqrt{8 \log n}}+O\left(\frac
{(\log\log n)^2}{(\log n)^{3/2}}\right),\\
\mathbb{V}\left(Z_{n-k} \right)&=\frac
{\pi^2/6-s_{k+1}^2}{2 \log n}+O\left(\frac
{1}{(\log n)^{2}}\right),
\end{align*}
where, for $k>0$, $s_k^1=\sum^k_{j=1} 1/j$, $s_k^2=\sum^k_{j=1} 1/j^2$ and $\gamma_0$ is the Euler constant.
\end{lemma}
\noindent \textbf{Proof of Lemma \ref{Cramer}.}
Following \cite{Cramer}, let $\xi_{n-k+1}=n(1-\Phi(Z_{n-k+1}))$ for $k\geqslant 1$. Since the random variables $\xi_{1}/n<...<\xi_{n}/n$ are the order statistics of $n$
independent uniform random variables, we see that $\xi_{n-k+1}$ has density
\[
f_{\xi_{n-k+1}}(x)=\left(
\begin{array}
[c]{c}%
n-1\\
k-1
\end{array}
\right)  \left(  \frac{x}{n}\right)  ^{k-1}\left(1-\frac{x}{n}\right)^{n-k}1_{\left[  0,n\right]  }(x).\]
\textbf{Step 1.} Write $\Gamma(k)=(k-1)!$ and observe that
\begin{align*}
\binom{n-1}{k-1}\left(  \frac{1}{n}\right)  ^{k-1} & =\exp\left(  \sum_{j=1}^{k-1}\log\left(  1-\frac{j}%
{n}\right)  \right)  \frac{1}{\Gamma(k)}
=\left(  1+O\left(  \frac{(\log n)^{3\theta}}{n}\right)  \right)  \frac
{1}{\Gamma(k)}
\end{align*}
since we have%
\begin{align*}
-\sum_{j=1}^{k}\frac{j}{n}-\sum_{j=1}^{k}\left(  \frac{j}{n}\right)  ^{2}  &
\leqslant\sum_{j=1}^{k}\log\left(  1-\frac{j}{n}\right)  \leqslant-\sum
_{j=1}^{k}\frac{j}{n}\\
\max_{1\leqslant k\leqslant C(\log n)^{\theta}}\left\vert \sum_{j=1}^{k}%
\log\left(  1-\frac{j}{n}\right)  +\sum_{j=1}^{k}\frac{j}{n}\right\vert  &
\leqslant\frac{1}{n}\sum_{j=1}^{\left[C  (\log
n)^{\theta}\right]}\frac{j^{2}}{n}=O\left(
\frac{(\log n)^{3\theta}}{n}\right)  .
\end{align*}
\noindent \textbf{Step 2.} For $k\geqslant 1$ we have

\begin{align*}
\mathbb{E}\left(  Z_{n-k+1}\right)    & =\mathbb{E}\left(  \Phi^{-1}\left(
1-\frac{\xi_{n-k+1}}{n}\right)  \right)  \\
& =\frac{(n-1)...(n-k+1)}{\Gamma(k)} \int_{0}^{n}
 \left(  \frac{x}{n}\right)  ^{k-1}\left(  1-\frac{x}{n}\right)
^{n-k}\Phi^{-1}\left(  1-\frac{x}{n}\right)  dx\\
& =
\left( 1-\frac{1}{n}\right)...\left( 1-\frac{k-1}{n}\right) \int_{0}^{n}\frac{x^{k-1}}{\Gamma(k)}\left(  1-\frac{x}%
{n}\right)  ^{n-k}\Phi^{-1}\left(  1-\frac{x}{n}\right)  dx\\
& =\exp\left(  -\frac{s_{k}^{1}}{n}-\frac{s_{k}^{2}}{2n^2}(1+o(1))\right) \left( E_{1,n}+E_{2,n} \right)
\end{align*}
where, for $p>\theta+1$, $x(n)=(\log n)^p$ and $\displaystyle f_{\Gamma
(k)}(x)=\frac{x^{k-1}}{\Gamma(k)}  e^{-x}$ for $x>0$, 
\begin{align*}
E_{1,n} & = (1+o(1)) \int_{0}^{x(n)}\Phi^{-1}\left(  1-\frac{x}{n}\right)  f_{\Gamma
(k)}(x)dx,\\
E_{2,n} & = \int_{x(n)}^{n}\Phi^{-1}\left(  1-\frac{x}{n}\right)\frac{x^{k-1}}{\Gamma(k)}  \left(  1-\frac{x}{n}\right)  ^{n-k}  dx.
\end{align*}
Assume that $k\leqslant C(\log n)^\theta$. By (\ref{phiinv}) it holds, for some $K>0$ and all n large enough,
\begin{align*}
|E_{2,n}| & \leqslant \int_{x(n)}^{n}x^{k-1} \left(  1-\frac{x}{n}\right)  ^{n-k} \left\vert \Phi^{-1}\left(  1-\frac{x}{n}\right) \right\vert  dx\\
& \leqslant K\sqrt{\log n}
\int_{x(n)}^{n/2}\exp \left( -(n-k)\frac{x}{n}+(k-1)\log x \right) dx\\
& + K \int_{n/2}^{n}x^{k-1} \left(  1-\frac{x}{n}\right)  ^{n-k} \sqrt{\log \left( \frac{1}{1-x/n} \right)}  dx\\
& \leqslant K \exp \left( -x(n) + \frac{C(\log n)^\theta + \log\log n}{2} + C(\log n)^{\theta+1} \right)\\
& + K \int_{n/2}^{n}x^{k-1} \left(  1-\frac{x}{n}\right)  ^{n-k-1}  dx\\
& \leqslant K \exp \left( -\frac{(\log n)^p}{2}\right) + K n^{k} \left(  \frac{1}{2}\right)  ^{n-k-1}\\
& \leqslant   K \exp \left( -(1+o(1))(\log n)^p\right).
\end{align*}
\noindent Now turn to%
\[
\int_{0}^{x(n)}\Phi^{-1}\left(  1-\frac{x}{n}\right)  \binom{n-1}{k-1}\left(
\frac{x}{n}\right)  ^{k-1}\left(  1-\frac{x}{n}\right)  ^{n-k}dx
\]
where, for $0<x<x(n)$, we have, by  (\ref{phiinv}),%
\begin{align}  
\Phi^{-1}(1-x/n)  & =  \sqrt{2\left(  \log\left(  n/x\right)  -\frac{1}{2}%
\log\log\left(  n/x\right)  -\frac{1}{2}\log(4\pi)-O\left(  \frac{\log
\log\left(  n/x(n)\right)  }{\log\left(  n/x(n)\right)  }\right)  \right)
}\nonumber\\
& =  \sqrt{2\log n}-\frac{2\log x+\log\log\left(  n/x\right)  +\log(4\pi
)}{2\sqrt{2\log n}}+O\left(  \frac{(\log\log n)^2}{(\log n)^{3/2}}\right)\label{phimoins1}
\end{align}
which is integrable near $0$ with respect to the above density since%
\begin{equation}\label{major}
0<\log(\log n-\log x)=\log\log n+\log\left(  1-\frac{\log x}{\log n}\right)
\leqslant\log\log n+\left\vert \frac{\log x}{\log n}\right\vert
\end{equation}
and $\log x,(\log x)^2$ are integrable with respect to any Gamma distribution. Hence%
\begin{eqnarray*}
E_{1,n} & = &\int_{0}%
^{x(n)}\left(  \sqrt{2\log n}-\frac{2\log x+\log\log n+\log(4\pi)}{\sqrt{8\log
n}}+O\left(  \frac{(\log\log n+\left\vert \log x\right\vert )^2}{(\log n)^{3/2}%
}\right)  \right)  f_{\Gamma(k)}(x)dx\\
& = &O\left(  \frac{(\log\log n)^2}{(\log n)^{3/2}}\right)  +\int_{0}^{x(n)}\left(
\sqrt{2\log n}-\frac{2\log x+\log\log n+\log(4\pi)}{\sqrt{8\log n}}\right)
f_{\Gamma(k)}(x)dx\\
& = &O\left(  \frac{(\log\log n)^2}{(\log n)^{3/2}}\right)  +\int_{0}^{+\infty
}\left(  \sqrt{2\log n}-\frac{2\log x+\log\log n+\log(4\pi)}{\sqrt{8\log n}%
}\right)  f_{\Gamma(k)}(x)dx
\end{eqnarray*}
since we have $x(n)=(\log n)^{p}$, $p>1$ thus, for any $s>1$,%
\[
\int_{x(n)}^{+\infty}f_{\Gamma
(k)}(x)dx=o\left(  \frac{1}{n^{s}}\right)
,\quad\int_{x(n)}^{+\infty}\log x f_{\Gamma(k)}(x)dx=o\left(  \frac{1}{n^{s}%
}\right)  .
\]
and moreover -- see \cite{Cramer} -- it holds
$$
\int_{0}^{+\infty}\log x\ f_{\Gamma(k)}(x)dx=s_{k+1}^{1}-\gamma_{0},
$$which yields the conclusion.\\
Similar computations give the claimed result for the variance. More precisely in the step 2 when substituing $\Phi^{-1}\left(  1-\frac{x}{n}\right)^2$ to $\Phi^{-1}\left(  1-\frac{x}{n}\right)$ in $E_{1,n}$ and $E_{2,n}$ it again appears that we can only consider integrals up to $x(n)$. Then it remains to compute, by substituing the expression of $\E (Z_{n-k})$ and using equation (\ref{phimoins1}) for $\Phi^{-1}\left(  1-\frac{x}{n}\right)$ :
\begin{eqnarray*}
& &\int_{0}^{x(n)}\left(\Phi^{-1}\left(  1-\frac{x}{n}\right)-\E (Z_{n-k})\right)^2  \binom{n-1}{k-1}\left(
\frac{x}{n}\right)  ^{k-1}\left(  1-\frac{x}{n}\right)  ^{n-k}dx\\
&= &\int_{0}^{x(n)}\left( -\frac{2(\log x-(s_{k+1}^1-\gamma_0))}{2\sqrt{2\log n}}+\frac{-\log\log\left(  n/x\right)  
+\log\log n}{2\sqrt{2\log n}}+O\left(\frac
{(\log\log n)^2}{(\log n)^{3/2}}\right)\right)^2\\
& \times&\binom{n-1}{k-1}\left(
\frac{x}{n}\right)  ^{k-1}\left(  1-\frac{x}{n}\right)  ^{n-k}dx.
\end{eqnarray*}
We conclude along the same lines as above by the upper bound (\ref{major}) and the fact that the variance of the logarithm of a variable with distribution $\Gamma(k)$ is $\displaystyle {\pi^2/6-s_{k+1}^2}$.
$\square$

\subsection{Proof of Theorem \ref{CVE}}

\noindent We intend to mimic the sheme of proof worked out in \cite{BFK17} and \cite{BF19} - specialized to the simpler case of the distance between the empirical and true c.d.f.'s instead of two correlated empirical ones. However all arguments have to be reconsidered since the almost sure controls by means of the law of the iterated logarithm and strong approximations can not be turned easily into $L_1$ controls. Indeed, what happens now is that the main part of the random integral we consider is also built from the extreme parts rather than the inner part only. Moreover, only a very short extreme interval can be neglected and the remainder extreme intervals define a divergent integral to be precisely evaluated as a series. This is why the expectation rate is no more a CLT rate. Note that the $\log\log n$ in this paper only comes from the primitive of $u(1-u)/h(u)^2$. Introduce the following decomposition, for $C>0$, $\gamma>1$ and $1<\theta\leqslant2$,
\begin{eqnarray*}
A_{n}&= &%
{\displaystyle\int\nolimits_{1-1/(n(\log n)^{\gamma})}^{1}}
\left(  Z_{n}-\Phi^{-1}(u)\right)  ^{2}du,\quad B_{n}=%
{\displaystyle\int\nolimits_{1-1/n}^{1-1/(n(\log n)^{\gamma})}}
\left(  Z_{n}-\Phi^{-1}(u)\right)  ^{2}du,\\
C_{n}&= &%
{\displaystyle\int\nolimits_{1-\left[  C(\log n)^{\theta}\right]/n}^{1-1/n}}
\left(  \mathbb{F}_{n}^{-1}(u)-\Phi^{-1}(u)\right)  ^{2}du,\quad D_{n}=
{\displaystyle\int\nolimits_{1/2}^{1-\left[  C(\log n)^{\theta}\right]/n}}
\left(  \mathbb{F}_{n}^{-1}(u)-\Phi^{-1}(u)\right)  ^{2}du.
\end{eqnarray*}
\textbf{Step 1.} We have, for $\gamma>1$,%
\[
\dfrac{nA_{n}}{\log\log n}\leqslant\dfrac{2Z_{n}^{2}}{(\log n)^{\gamma}%
\log\log n}+\dfrac{2n}{\log\log n}%
{\displaystyle\int\nolimits_{1-1/(n(\log n)^{\gamma})}^{1}}
\left(  \Phi^{-1}(u)\right)  ^{2}du
\]
where%
\[
\lim_{n\rightarrow+\infty}\dfrac{\mathbb{E}\left(  Z_{n}^{2}\right)  }{(\log
n)^{\gamma}\log\log n}=0
\]
and%
\begin{eqnarray*}%
{\displaystyle\int\nolimits_{1-1/(n(\log n)^{\gamma})}^{1}}
\left(  \Phi^{-1}(u)\right)  ^{2}du &  = &%
{\displaystyle\int\nolimits_{1-1/(n(\log n)^{\gamma})}^{1}}
2\log\left(  \frac{1}{1-u}\right)  \left(  1+o(1-u)\right)  ^{2}du\\
&  = &\left[  -2(1-u)\log\left(  \frac{1}{1-u}\right)  \right]  _{1-1/(n(\log
n)^{\gamma})}^{1}=O\left(  \frac{1}{n(\log n)^{\gamma-1}}\right)
\end{eqnarray*}
hence
\[
\lim_{n\rightarrow+\infty}\dfrac{n\mathbb{E}\left(  A_{n}\right)  }{\log\log
n}=0.
\]
\textbf{Step 2.} Notice that for all $u\in [1-1/n,1-1/(n(\log n)^\gamma)]$, we have 
$$\Phi^{-1}(u)=\sqrt{2\log n}+O\left(\frac{\log\log n}{\sqrt{\log n}}\right).$$
Next observe that
\begin{align*}
\mathbb{E}\left(  B_{n}\right)   &  =\frac{\mathbb{V}\left(  Z_{n}\right)
}{n}\left(  1-\frac{1}{(\log n)^{\gamma}}\right)  +%
{\displaystyle\int\nolimits_{1-1/n}^{1-1/(n(\log n)^{\gamma})}}
\left(  \mathbb{E}\left(  Z_{n}\right)  -\Phi^{-1}(u)\right)  ^{2}du\\
&  =O\left(  \frac{1}{n\log n}\right)  +O\left(  \frac{(\log\log n)^{2}}{n\log
n}\right),
\end{align*}
hence%
\[
\lim_{n\rightarrow+\infty}\dfrac{n\mathbb{E}\left(  B_{n}\right)  }{\log\log
n}=0.
\]
\textbf{Step 3.} Start with%
\[
C_{n}=\sum_{k=1}^{\left[  C(\log n)^{\theta}\right]  }%
{\displaystyle\int\nolimits_{1-(k+1)/n}^{1-k/n}}
\left(  Z_{n-k}-\Phi^{-1}(u)\right)  ^{2}du.
\]
Recall that%
\[
s_{k}^{1}-\gamma_{0}=\log k+\frac{1}{2k}+O\left(  \frac{1}{k^{2}}\right)  .
\]
Now, for $1\leqslant k\leqslant\left[  C(\log n)^{\theta}\right]  $ and
$u\in\left[  1-(k+1)/n,1-k/n\right]  $ we have%
\begin{align*}
\Phi^{-1}(u)  & =\sqrt{2\left(  \log\left(  1-u\right)  -\frac{1}{2}\log
\log\left(  1-u\right)  -\frac{1}{2}\log(4\pi)-O\left(  \frac{\log\log n}{\log
n}\right)  \right)  }\\
& =\sqrt{2\log n}-\frac{2\log k+\log\log n+\log(4\pi)}{\sqrt{8\log n}%
}+O\left(  \frac{(\log\log n)^2}{(\log n)^{3/2}}\right)
\end{align*}
thus, by Lemma \ref{Cramer}, we have, uniformly in $k$,%
\[
\mathbb{V(}Z_{n-k})=\frac{\pi^{2}/6-s_{{k+1}}^{2}}{2\log n}+O\left(  
\frac{1}{(\log n)^2} \right)
\]
then
\begin{align*}
& \mathbb \E\left(  \left(  Z_{n-k}-\Phi^{-1}(u)\right)  ^{2}\right)\\
&= \mathbb V(Z_{n-k})+\left(  \mathbb{E(}Z_{n-k})-\Phi^{-1}(u)\right)  ^{2}\\
&= \frac{\pi^{2}/6-s_{k+1}^{2}}{2\log n}+O\left( \frac{1}{(\log n)^{2}}\right)  +\left(  \frac{\log k-(s_{k+1}^{1}-\gamma_{0})}%
{\sqrt{2\log n}}+O\left(  \frac{(\log\log n)^2}{(\log n)^{3/2}}\right)  \right)  ^{2}\\
&= \frac{\pi^{2}/6-s_{k+1}^{2}}{2\log n}+O\left(   \frac{1}{(\log n)^{2}}  \right)  +\left(  \frac{1+O(1/k)}{2k\sqrt{2\log n}}+O\left(
\frac{(\log\log n)^2}{(\log n)^{3/2}}\right)  \right)  ^{2}.
\end{align*}
As a consequence,
\begin{align*}
& \mathbb{E}\left(  C_{n}\right)\\
& =\frac{1}{n}\sum_{k=1}^{\left[  C(\log
n)^{\theta}\right]  }\frac{\pi^{2}/6-s_{k+1}^{2}}{2\log n}+O\left(  \frac
{1}{n(\log n)^{2-\theta}}\right)  +\frac{1}{n}\sum
_{k=1}^{\left[  C(\log n)^{\theta}\right]  }\left(  \frac{1+O(1/k)}%
{2k\sqrt{2\log n}}+O\left(  \frac{(\log\log n)^2}{(\log n)^{3/2}}\right)  \right)  ^{2}\\
& = O\left(  \frac{(\log n)^{\theta/2}}{n\log n}\right)  +\frac{1}{n}%
\sum_{k=\left[  (\log n)^{\theta/2}\right]  }^{\left[  C(\log n)^{\theta
}\right]  }\frac{\pi^{2}/6-s_{k+1}^{2}}{\log n}+O\left(  \frac{1}{n\log
n}\right)  +O\left(  \frac{(\log\log n)^{3}}{n(\log n)^{3-\theta}}\right)  \\
& \leqslant \frac{C(\log n)^{\theta}}{n\log n}\sum_{j=\left[  (\log n)^{\theta
/2}\right]  }^{+\infty}\frac{1}{j^{2}}+O\left(  \frac{(\log n)^{\theta/2}%
}{n\log n}\right)  \\
& = O\left(  \frac{(\log n)^{\theta/2}}{n\log n}\right).
\end{align*}
Thus, for any $\theta\leqslant 2$ we have 
\[
\lim_{n\to +\infty}\frac{n \mathbb{E} (C_n)}{\log\log n}=0.
\]

\noindent\textbf{Step 4.} Now we compute the limit of the main deterministic contribution to the main stochastic term $D_n$, namely
\[
D_{1,n}=
{\displaystyle\int\nolimits_{1/2}^{1-\left[  C(\log n)^{\theta
}\right]/n}}
\frac{u(1-u)}{h^{2}(u)}du.
\]
Let $v_{n}$ be such that $\log v_{n}=(\log n)^{\varepsilon_{n}}$, $\displaystyle \lim_{n\to +\infty}\varepsilon_n=0,
 \lim_{n\to +\infty}\varepsilon_{n} \log\log
n=+\infty$. By using (\ref{h}) it holds
\begin{align*}
& \dfrac{1}{\log\log n}%
{\displaystyle\int\nolimits_{1-1/v_{n}}^{1-\left[  C(\log n)^{\theta
}\right]/n}}
\frac{u(1-u)}{h^{2}(u)}du  \\
& = \dfrac{1+o(1)}{2\log\log n}(\log(\log n-\log(\left[  C(\log n)^{\theta
}\right]))-\log\log v_{n})\\
&  = \dfrac{1+o(1)}{2\log\log n}\log\left(  \frac{(1+o(1))\log n}{\log v_{n}}\right) \\
& = \dfrac{1+o(1)}{2}(1-\varepsilon_{n})
\end{align*}
and %
\[
\dfrac{1}{\log\log n}%
{\displaystyle\int\nolimits_{1/2}^{1-1/v_{n}}}
\frac{u(1-u)}{h^{2}(u)}du\leqslant\dfrac{1+o(1)}{2\log\log n}%
(\log\log v_{n})=\dfrac{1+o(1)}{2}\varepsilon_{n}.\] 

Therefore \begin{equation}\label{8pi2}
\lim_{n\to +\infty} \dfrac{D_{1,n}}{\log\log n}=\dfrac{1}{2}.
\end{equation}
Compared with the result of \cite{Bickel78} recalled at Remark \ref{Rem0} the truncation at level $1/v_n$ instead of $1/n$ preserves the same first order.\medskip

\noindent\textbf{Step 5.} To show that $\E(D_n)$ behaves as $D_{1,n}+o(1)$ we proceed as in \cite{BFK17} with strong approximation arguments. First, we substitute the uniform quantile process to the general quantile process with a sharp control of the expectation of the random error terms in the Taylor Lagrange expansion. For short, write $d_n=\left[  C(\log n)^{\theta
}\right]/n$ and $\beta
_{n}^{X}(u)=\sqrt{n}( \mathbb{F}_{n}^{-1}(u)-\Phi^{-1}(u))$ so that
\[
\frac{nD_{n}}{\log\log n}=\frac{1}{\log\log n}\int_{1/2}^{1-d_{n}}(\beta
_{n}^{X}(u))^{2}du.
\]
Defining $U_{i}=\Phi(X_{i})$ which is uniform on $(0,1)$ we obviously have $U_{(i)}=\Phi(X_{(i)})$. Let denote $\mathbb{F}_{n}^{U}$ the uniform empirical $c.d.f.$ associated to the $U_i$ and define the underlying uniform quantile process to be
\[
\beta_{n}^{U}(u)=\sqrt{n}((\mathbb{F}_{n}^{U})^{-1}(u)-u)=\sqrt{n}%
(\Phi(\mathbb{F}_{n}^{-1}(u))-u).
\]
Thus for all $1/2\leqslant u\leqslant 1-d_{n}$ there exists a random $u^{\ast}$ such that $\left\vert
u-u^{\ast}\right\vert \leqslant\left\vert \beta_{n}^{U}(u)\right\vert
/\sqrt{n}$ and
\begin{align*}
\beta_{n}^{X}(u)h(u) &  =\sqrt{n}(\mathbb{F}_{n}^{-1}(u)-\Phi^{-1}(u))h(u)\\
&  =\sqrt{n}(\Phi^{-1}(\Phi(\mathbb{F}_{n}^{-1}(u)))-\Phi^{-1}(u))h(u)\\
&  =\sqrt{n}\left(  \frac{\Phi(\mathbb{F}_{n}^{-1}(u))-u}{h(u)}+\frac
{h^{\prime}(u^{\ast})}{2h^{2}(u^{\ast})}(\Phi(\mathbb{F}_{n}^{-1}%
(u))-u)^{2}\right)  h(u)\\
&  =\beta_{n}^{U}(u)+r_{n}(u)
\end{align*}
with
\begin{align*}
r_{n}(u) &  =\frac{1}{2\sqrt{n}}\left(  \beta_{n}^{U}(u)\right)  ^{2}%
\frac{h^{\prime}(u^{\ast})}{h(u^{\ast})}\frac{h(u)}{h(u^{\ast})}\\
&  =\frac{1}{2\sqrt{n}}\left(  \frac{\beta_{n}^{U}(u)}{\sqrt{1-u}}\right)
^{2}\left(  \frac{1-u}{1-u^{\ast}}\right)  \left(  (1-u^{\ast})\frac
{\Phi^{\prime\prime}(\Phi^{-1}(u^{\ast}))}{{\Phi^{\prime}}^{2}(\Phi^{-1}(u^{\ast}))}\right)
\frac{h(u)}{h(u^{\ast})}.
\end{align*}
We study
\[
\frac{nD_{n}}{\log\log n}=\frac{1}{\log\log n}\int
_{1/2}^{1-d_{n}}(\beta_{n}^{U}(u)+r_{n}(u))^{2}\frac{du}{h(u)^{2}}.
\]
Since we have%
\begin{align*}
\sup_{0<u<1}u(1-u)\frac{\left\vert \Phi^{\prime\prime}(\Phi
^{-1}(u))\right\vert }{{\Phi^{\prime}}^{2}(\Phi^{-1}(u))}=1
\end{align*}
it holds, by Lemma 6.1.1 in \cite{CH93},%
\[
0\leqslant\frac{h(u)}{h(u^{\ast})}\leqslant\frac{\max(u,u^{\ast})}%
{\min(u,u^{\ast})}\frac{1-\min(u,u^{\ast})}{1-\max(u,u^{\ast})}.
\]
Now we introduce the sequence of events, with $0<\varepsilon <1$,%
\begin{equation}\label{An}
{\cal A}_{n}=\left\{  \left\vert \frac{\beta_{n}^{U}(u)}{\sqrt{u(1-u)}%
}\right\vert \leqslant(1-\varepsilon)\sqrt{n(1-u)},\ d_n<u<1-d_{n}\right\}.
\end{equation}
\noindent On the event ${\cal A}_n$ we have the following control of $u^{\ast}$,
\[
\frac{\max(u,u^{\ast})}%
{\min(u,u^{\ast})}\frac{1-\min(u,u^{\ast})}{1-\max(u,u^{\ast})}\leqslant \frac{4}{\varepsilon^2}
\]
since, for instance,
\[
0\leqslant\frac{1-u}{1-u^{\ast}}\leqslant1+\frac{u^{\ast}-u}{1-u-(u^{\ast}-u)}        \leqslant 1+\frac{\left\vert \frac{\beta_{n}^{U}%
(u)}{\sqrt{u(1-u)}}\frac{1}{\sqrt{n(1-u)}}\right\vert }{1-\left\vert
\frac{\beta_{n}^{U}(u)}{\sqrt{u(1-u)}}\frac{1}{\sqrt{n(1-u)}}\right\vert
}\leqslant \frac{2}{\varepsilon},
\]
\[
0 \leqslant\frac{u}{u^{\ast}}=1+\frac{u-u^{\ast}}{u+u^{\ast}-u}\leqslant
1+\frac{\left\vert \frac{\beta_{n}^{U}%
(u)}{\sqrt{u(1-u)}}\frac{1}{\sqrt{n(1-u)}}\right\vert }{1-\left\vert
\frac{\beta_{n}^{U}(u)}{\sqrt{u(1-u)}}\frac{1}{\sqrt{n(1-u)}}\right\vert
}\leqslant \frac{2}{\varepsilon},
\]
and the same holds for the reverse ratios. Hence we have
\[
1_{{\cal A}_n} r_{n}(u) \leqslant\frac{4}{\varepsilon ^3\sqrt{n}}\left(  \frac{\beta_{n}^{U}(u)}{\sqrt{1-u}}\right)
^{2}
\]
thus
\begin{align*}
\E \left(\int_{1/2}^{1-d_{n}}1_{{\cal A}_n}\frac{r_{n}(u)^{2}}{h(u)^{2}} du\right) &\leqslant \int_{1/2}^{1-d_{n}}\frac{16}{\varepsilon ^6 n (1-u)}
\E\left(  \frac{\beta_{n}^{U}(u)}{\sqrt{1-u}}\right)
^{4} \frac{1-u}{h(u)^{2}}du.\end{align*}
By Lemma \ref{moments} below and (\ref{8pi2}) we have, when $\theta=2$, \begin{equation}\label{step4}
\sup_{1/2<u<1-d_n}\E\left(  \frac{\beta_{n}^{U}(u)}{\sqrt{1-u}}\right)
^{4}=O(1), \quad \int_{1/2}^{1-d_{n}} \frac{1-u}{h(u)^{2}}du=O(\log\log n).\end{equation} It ensues
\begin{align*}
\E \left(\int_{1/2}^{1-d_{n}}1_{{\cal A}_n}\frac{r_{n}(u)^{2}}{h(u)^{2}} du\right)=O\left(\frac{\log\log n}{(\log n)^2}\right).
\end{align*}
By using the Cauchy-Schwartz inequality
 we easily get \begin{equation*}
      \displaystyle\lim_{n\to+\infty}\E\left(\int_{1/2}^{1-d_{n}}1_{{\cal A}_n}\frac{\beta_{n}^{U}(u)r_{n}(u)}{h(u)^{2}}du\right)=~0, 
 \end{equation*}
 since by (\ref{step4}) we have, again for $\theta=2$,
     $$\displaystyle \int_{1/2}^{1-d_{n}} \frac{\E(\beta_{n}^{U}(u)^2)}{h(u)^{2}}du=O(\log\log n).$$
     
\medskip     
\noindent\textbf{Step 6.} Next we evaluate the probability of the rare event ${\cal A}_n^c$ from (\ref{An}). To this aim we work on the KMT probability space where we can define a sequence $\mathbb B_n$ of standard Brownian bridges approximating the processes $\beta_n^U$ in such a way that the error process $w_n=\beta_n^U-\mathbb B_n$ satisfies, for universal positive constants $c_1, c_2, c_3$ and all $x>0, n \geqslant 1$,
\begin{equation}\label{KMT}
\mathbb{P}\left(  \sup_{0<u<1}\left\vert w_{n}(u)\right\vert
>\frac{c_{1}}{\sqrt{n}}\left(x+\log n\right)  \right)\leqslant c_{2}\exp(-c_{3}x).   
\end{equation}
Hence we have
\begin{align*}
\mathbb P({\cal A}_n^c) &= \mathbb{P}\left( \exists u\in[1/2,1-d_{n}], \left\vert \frac{\beta_{n}^{U}(u)}{\sqrt{u(1-u)}%
}\right\vert >(1-\varepsilon)\sqrt{n(1-u)}%
\right)  \\
& \leqslant\mathbb{P}\left(  \sup_{1/2<u<1-d_{n}}\left\vert \frac{\beta
_{n}^{U}(u)}{\sqrt{u(1-u)}}\right\vert >(1-\varepsilon)(\log n)^{\theta
/2}\right)  \\
& \leqslant\mathbb{P}\left( \left\{  \sup_{1/2<u<1-d_{n}}\left\vert
\frac{\mathbb B_{n}(u)}{\sqrt{u(1-u)}}\right\vert >\frac{1-\varepsilon}{2}(\log
n)^{\theta/2}\right\} \dots\right.  \\
& \left.  \dots\cap\left\{  \sup_{1/2<u<1-d_{n}}\left\vert \frac
{w_{n}(u)}{\sqrt{u(1-u)}}\right\vert \leqslant
\frac{1-\varepsilon}{2}(\log n)^{\theta/2}\right\} \right) \\
& +\mathbb{P}\left(  \sup_{1/2<u<1-d_{n}}\left\vert \frac{w_{n}(u)}{\sqrt{u(1-u)}}\right\vert >\frac{(1-\varepsilon)}%
{2}(\log n)^{\theta/2}\right)  \\
& \leqslant\mathbb{P}\left(  \sup_{1/2<u<1-d_{n}}\left\vert \frac{\mathbb B_{n}%
(u)}{\sqrt{u(1-u)}}\right\vert >\frac{1-\varepsilon}{2}(\log n)^{\theta
/2}\right)\\
& +\mathbb{P}\left(  \sup_{1/2<u<1-d_{n}}\left\vert w_{n}
(u)\right\vert >\sqrt{C}\frac{1-\varepsilon}%
{2}\frac{(\log n)^{\theta}}{\sqrt{n}}\right).
\end{align*}
Recall that $1<\theta \leqslant 2$. By the theorem of Borell-Sudakov (see \cite{Bor75}, \cite{ledtal}) and (\ref{KMT}) we obtain, for any $\gamma>2$,  the constant $C$ fixed as large as needed and all $n$ large enough,
\begin{align*}
\mathbb P({\cal A}_n^c) & \leqslant\exp\left(  -\frac{(1-\varepsilon)^{2}(\log n)^{\theta}}%
{8\sup_{1/2<u<1-d_{n}}(\Var(B_{n}(u)/\sqrt{u(1-u)})})\right)
+c_2\exp\left(  -c_3(\log n)^{\theta}\right) \\
& \leqslant\exp\left(  -\frac{(1-\varepsilon)^{2}}{8}(\log n)^{\theta}\right)
+c_2\exp\left(  -c_3(\log n)^{\theta}\right)\\
& \leqslant\frac{1}{n^{\gamma}}.
\end{align*}
Therefore we get, for any $0<b<\gamma / 2-1$,
\begin{align*}
&  \mathbb{E}\left(  1_{{\cal A}_n^c}\int_{1/2}^{1-d_{n}}n(F_{n}^{-1}(u)-\Phi
^{-1}(u))^{2}du\right)  \\
&  \leqslant\mathbb{P}\left(  {\cal A}_n^c\right)  2n\int_{0}^{1}\Phi^{-1}%
(u)^{2}du+2\mathbb{E}\left(  1_{{\cal A}_n^c}nZ_{n}^{2}\right)  \\
&  \leqslant2n\mathbb{P}\left(  {\cal A}_n^c\right)  +\sqrt{\mathbb{P}\left(
A_{n}^c\right)  n^{2}\mathbb{E}\left(  Z_{n}^{4}\right)  }=O\left(\frac{1}{n^b}\right).
\end{align*}

\noindent\textbf{Step 7.} It remains to study $$\displaystyle \frac{1}{\log\log n}\int
_{1/2}^{1-d_{n}}\E(\beta_{n}^{U})^{2}\frac{du}{h(u)^{2}}.$$
At this stage the approximation bounds play a crucial role and there is no room for relaxing the trimming constraints. To be more specific the only allowed choice $\theta \leqslant 2$ is $\theta=2$. Choose an arbitrarily large constant $C>0$. Given any $0<\eta<1$, consider the sequence of events
\[
\mathcal{B}_{n}=\left\{  \left\vert w_{n}(u)\right\vert <\eta
\sqrt{u(1-u)}, \frac{1}{2}<u<1-d_{n}\right\}.
\]
By (\ref{KMT}),
for any $k_{1}>0$ there exists $C=C_{\eta
}>(1+k_{1}/c_{3})^{2}/\eta^{2}>0$ and $n_{0}>0$ large enough
such that for all $n>n_{0}$ we have%
\begin{align*}
1-\mathbb{P}\left(  \mathcal{B}_{n}\right)    & \leqslant\mathbb{P}\left(
\sup_{1/2<u<1-d_{n}}\left\vert w_{n}(u)\right\vert >\eta\sqrt
{\frac{C_{\eta}}{n}}(\log n)^{\theta/2}\right)  \\
& \leqslant\mathbb{P}\left(  \sup_{0<u<1}\left\vert w_{n}(u)\right\vert
>\frac{c_{1}}{\sqrt{n}}\left(  (\eta\sqrt{C_{\eta}}-1)\log
n+\log n\right)  \right)  \\
& \leqslant c_{2}\exp(-c_{3}(\eta\sqrt{C_{\eta}}-1)\log n)\\
& \leqslant\frac{1}{n^{k_{1}}}.
\end{align*}
\begin{lemma}\label{moments} For any $p\geqslant 1$ there exist constants $C>0$ and $\kappa_p$ such that we have, for $d_n=[C\frac{(\log n)^2}{n}]$ and all $n$ large enough,
\[\sup_{d_n<u<1-d_{n}}\E\left( \frac{|w_{n}(u)|}{\sqrt{u(1-u)}}\right)^{p}<2\eta^p,\quad \sup_{d_n<u<1-d_{n}}\E\left( \frac{|\beta_{n}^{U}(u)|}{\sqrt{u(1-u)}}\right)^{p}<\kappa_p. \]

\end{lemma}

\noindent \textbf{Proof of Lemma \ref{moments}.} Start with%

\[
\mathbb{E}\left(  \frac{\left\vert w_{n}(u)\right\vert}{\sqrt{u(1-u)}}\right)  ^{p}\leqslant
\mathbb{\eta}^{p}+\mathbb{E}\left(  1_{\mathcal{B}_{n}^{c}}\frac{\left\vert w_{n}(u)\right\vert}{\sqrt{u(1-u)}}\right)  ^{p}
\]
then set, for $k\geqslant 0$,%
\begin{align*}
\mathcal{F}_{n}  & =\left\{  \left\vert \frac{B_{n}(u)}{\sqrt{u(1-u)}%
}\right\vert <n:d_{n}<u<1-d_{n}\right\}  ,\\
\mathcal{F}_{n}^{c}  & \subset%
{\displaystyle\bigcup\nolimits_{k\in\mathbb{N}}}
\mathcal{F}_{n,k},\\
\mathcal{F}_{n,k}  & =\left\{  n+k\leqslant\sup_{0<u<1}\left\vert \frac
{B_{n}(u)}{\sqrt{u(1-u)}}\right\vert <n+k+1\right\}  .
\end{align*}
Since $\left\vert \beta_{n}^{U}(u)/\sqrt{u(1-u)}\right\vert \leqslant n$ for
$d_{n}<u<1-d_{n}$ and all $n$ large enough, we have%
\[
1_{\mathcal{F}_{n,k}}\sup_{d_{n}<u<1-d_{n}}\frac{\left\vert w_{n}%
(u)\right\vert }{\sqrt{u(1-u)}}\leqslant2n+k+1,\quad1_{\mathcal{F}_{n}}%
\sup_{d_{n}<u<1-d_{n}}\frac{\left\vert w_{n}(u)\right\vert }{\sqrt{u(1-u)}%
}\leqslant2n.
\]
By Sudakov-Borell theorem it holds $\mathbb{P}\left(  \mathcal{F}%
_{n,k}\right)  \leqslant\exp\left(  -(n+k)^{2}/2\right)  $ whereas
$\mathbb{P}\left(  \mathcal{B}_{n}^{c}\right)  <1/n^{k_{1}}$. Hence by choosing $k_1>p$ it holds
\begin{align*}
\mathbb{E}\left(  1_{\mathcal{B}_{n}^{c}}\frac{\left\vert w_{n}(u)\right\vert}{\sqrt{u(1-u)}%
}\right)  ^{p}  & \leqslant\left(
{\displaystyle\sum\limits_{k\in\mathbb{N}}}
(2n+k+1)^{p}\mathbb{P}\left(  \mathcal{F}_{n,k}\right)  \right)
+\mathbb{E}\left(  (2n)^{p}1_{\mathcal{F}_{n}\cap\mathcal{B}_{n}^{c}}\right)
\\
& \leqslant%
{\displaystyle\sum\limits_{k\in\mathbb{N}}}
(2n+k+1)^{p}\exp\left(  -(n+k)^{2}/2\right)  +\frac{(2n)^{p}}{n^{k_{1}}}\\
& =o(1), 
\end{align*}
which proves the first claimed upper bound. Since $$\E\left( |B_{n}(u)|/\sqrt{u(1-u)}\right)^{p}<+\infty$$ doesn't depend on $n$ the second expectation bound follows.$\quad \square$

\noindent By Lemma \ref{moments} we get%
\begin{align*}
&\frac{1}{\log\log n}\mathbb{E}\left(  \int_{1/2}^{1-d_{n}}%
(w_{n}(u))^{2}\frac{du}{h(u)^{2}}\right) \\
&=\frac{1}{\log\log n}\mathbb{E}%
\left(  \int_{1/2}^{1-d_{n}}\left(  \frac{w_{n}(u)}{\sqrt{u(1-u)}%
}\right)  ^{2}\frac{u(1-u)}{h(u)^{2}}du\right)  =O(\eta^{2})
\end{align*}
\noindent and, by (\ref{8pi2}),
\begin{align*}
& \frac{1}{\log\log n}\mathbb{E}\left(  \left\vert \int_{1/2}^{1-d_{n}}%
\frac{w_{n}(u) \mathbb B_{n}(u)}{h(u)^{2}}du\right\vert \right)  \\
& \leqslant\sqrt{\frac{1}{\log\log n}\mathbb{E}\left(  \int_{1/2}^{1-d_{n}%
}\frac{w_{n}(u)^{2}}{h(u)^{2}}du\right)  }\sqrt{\frac{1}{\log\log n}%
\mathbb{E}\left(  \int_{1/2}^{1-d_{n}}\frac{(\mathbb B_{n}(u))^{2}}{h(u)^{2}%
}du\right)  }\\%
& =O(\eta)\sqrt{\frac{1}{\log\log n}\int_{1/2}^{1-d_{n}%
}\frac{u(1-u)}{h(u)^{2}}du}.%
\end{align*}
By choosing $\eta$ as small as desired, the first assertion of Theorem \ref{CVE} is proved. \bigskip

\noindent\textbf{Step 8.} 
 The sequence $\displaystyle \sqrt{n / \log\log n}
W_{2}(\F_{n},\Phi)$ is bounded in $L^2$, thus uniformly integrable, and from  (\ref{Del}) (see \cite{DelBarrio05}) converges in probability to $1$. Thus the convergence holds in $L^1$, which establishes the second assertion of Theorem \ref{CVE}.$\quad \square$ \bigskip

\subsection{Proof of Theorem \ref{Theo2}}

\noindent In Theorem 11 of \cite{BF19} we proved that  $\displaystyle nW_2^2(\F_n,\G_n)$ converges in distribution to $$\|{\cal G}\|_2^2 = \displaystyle \int_0^1 \left(\frac{\mathbb{B}^{X}(u)}{h(u)}-\frac{\mathbb{B}^{Y}(u)}{h(u)}\right)^2 du$$
under assumptions on the common probability distribution $F$ of the samples ensuring that $\sqrt{n}(\F_n^{-1}(u)-F^{-1}(u))$ and $\sqrt{n}(\G_n^{-1}(u)-G^{-1}(u))$ can be simultaneously approximated on a suitable sub-interval of $[0,1]$ by $\mathbb{B}^{X}(u)/h(u)$ and $\mathbb{B}^{Y}(u)/h(u)$ respectively. Here $ \mathbb{B}^X(u)$ and $\mathbb{B}^Y(u)$ are two standard Brownian bridges coupled to the marginal samples respectively, and are then correlated together as mentionned at Theorem \ref{Theo2} if the two samples are. In \cite{BF19} the imposed assumptions for the Gaussian approximation concerned the tail of $F$ with respect to the cost function, and the integrability condition
$$\int_0^1\frac{u(1-u)}{h^2(u)}du<+\infty$$ was morerover required.
Under the latter condition, the expectation of $\|{\cal G}\|_2^2$ is finite since it is bounded by $ 4\int_0^1 u(1-u)/h^2(u)du$. Now, this upper bound is appropriated to the independent case whereas in our currently dependent case the sample is Gaussian and $$\mathbb{E} (\|{\cal G}\|_2^2)=2 \int_0^1\frac{u-C_\rho(u)}{h^2(u)}du$$which we shall next prove to be finite if $0<|\rho|<1$. Then, as the tail conditions of Theorem 11 in \cite{BF19} are satisfied by the Gaussian distribution $F=G=\Phi$, the weak convergence of $n W_{2}^2(\mathbb F_n,\mathbb G_n)$ is easily established by a straightforward adaptation of the proof of the latter theorem. This long and technical proof is thus omitted. Notice that in the case $\rho=0$ we have $\mathbb{E}(\|{\cal G}\|_2^2) = 2\int_0^1u(1-u)/h^2(u)du=+\infty$ and therefore by \cite{CHS93} the random variable $\|{\cal G}\|_2^2=+\infty\ a.s.$ and $nW_2^2(\F_n,\G_n)$ do not weakly converges.\\

Let us prove that
$$\displaystyle \int_0^1\frac{u-C_\rho(u)}{h^2(u)}du<+\infty.$$
Notice that for $a>0$, as $u\to 1$, 
$$ 1-\Phi(a\Phi^{-1}(u))=(4\pi)^{\frac{1-a^2}{2}}\frac{(1-u)^{a^2}}{ a(\log(\frac{1}{1-u}))^{\frac{1-a^2}{2}}}\left(1+O\left(\frac{\log\log(\frac{1}{1-u})}{\log (\frac{1}{1-u})}\right) \right).$$

\noindent First assume that $-1<\rho<0$. It holds

\begin{align*}
u-C_\rho(u)&=u-\frac{1}{2\pi \sqrt{1-\rho ^{2}}}\int_{-\infty
}^{\Phi ^{-1}(u)}\int_{-\infty }^{\Phi ^{-1}(u)}\exp \left( -\frac{%
x^{2}+y^{2}-2\rho xy}{2(1-\rho ^{2})}\right) dxdy\\
&=u-\int_{-\infty
}^{\Phi ^{-1}(u)}\frac{1}{\sqrt{2\pi}}e^{-\frac{y^2}{2}}\Phi\left(\frac{\Phi^{-1}(u)-\rho y}{\sqrt{1-\rho ^{2}}}\right)dy\\
&=u-\int_0^u\Phi\left(\frac{\Phi^{-1}(u)-\rho\Phi^{-1}(v)}{\sqrt{1-\rho ^{2}}}\right)dv=\int_0^u \left(1-\Phi\left(\frac{\Phi^{-1}(u)-\rho\Phi^{-1}(v)}{\sqrt{1-\rho ^{2}}}\right)\right)dv\\
&=\int_0^u \left(1-\Phi\left(\Phi^{-1}(u)\sqrt{\frac{1-\rho}{1+\rho}}+\frac{\rho(\Phi^{-1}(u)-\Phi^{-1}(v))}{\sqrt{1-\rho ^{2}}}\right)\right)dv\\
&\leqslant u\left(1-\Phi\left(\Phi^{-1}(u)\sqrt{\frac{1-\rho}{1+\rho}}\right)\right)
=O\left( \frac{u(1-u)^{\frac{1-\rho}{1+\rho}}}{(\log(\frac{1}{1-u}))^{\frac{2\rho}{1+\rho}}}\right), \quad u\to 1,
\end{align*}
which proves that $(u-C_\rho(u))/h^2(u)$ is integrable near $1$ since $-1<\rho<0$. By symmetry the same holds near $0$.\\ 
Next the case $0<\rho<1$ near $1$ follows from the equality
\begin{align*}
u-C_\rho(u)&=\int_0^u \left(1-\Phi\left(\frac{\Phi^{-1}(u)-\rho\Phi^{-1}(v)}{\sqrt{1-\rho ^{2}}}\right)\right)dv\\
&=\int_0^{\frac{1}{2}} \left(1-\Phi\left(\frac{\Phi^{-1}(u)-\rho\Phi^{-1}(v)}{\sqrt{1-\rho ^{2}}}\right)\right)dv+
\int_{\frac{1}{2}}^{u} \left(1-\Phi\left(\frac{\Phi^{-1}(u)-\rho\Phi^{-1}(v)}{\sqrt{1-\rho ^{2}}}\right)\right)dv.\\
\end{align*}
Then we get, for the first term, the upper bound
\begin{align*}
 \int_{0}^{\frac{1}{2}} \left(1-\Phi\left(\frac{\Phi^{-1}(u)-\rho\Phi^{-1}(v)}{\sqrt{1-\rho ^{2}}}\right)\right)dv\leqslant \frac{1}{2} \left(1-\Phi\left(\frac{\Phi^{-1}(u)}{\sqrt{1-\rho ^{2}}}\right)\right)   
\end{align*}
that is, up to a logarithmic factor, of order 
$(1-u)^{\frac{1}{1-\rho^2}}$ as $u\to 1$. 

\noindent The second term needs more attention. First we choose $0<\alpha<1$ such that for all $v\in[1/2,1-(1-u)^{\alpha^2}]$ we have, for $u$ close to $1$ and $\eta$ arbitrarily small, $\Phi^{-1}(v)\leqslant(\alpha+\eta)\Phi^{-1}(u)$ and $1-\alpha\rho>\sqrt{1-\rho ^{2}}$. We take $\alpha<(1-\sqrt{1-\rho ^{2}})/\rho$, which is actually less than $\rho$ and we have for $u$ close enough to $1$,
$$\Phi^{-1}(v)\leqslant\Phi^{-1}(1-(1-u)^{\alpha^2})\leqslant (\alpha+\eta)\Phi^{-1}(u).$$
Thus it comes
\begin{align*}
\int_{\frac{1}{2}}^{1-(1-u)^{\alpha^2}}
\left(1-\Phi\left(\frac{\Phi^{-1}(u)-\rho\Phi^{-1}(v)}{\sqrt{1-\rho ^{2}}}\right)\right)dv\leqslant \frac{1}{2}\left(1-\Phi\left(\frac{(1-(\alpha+\eta)\rho)\Phi^{-1}(u)}{\sqrt{1-\rho ^{2}}}\right)\right)
\end{align*}
that is, up to a logarithmic factor, of order $(1-u)^{\frac{(1-(\alpha+\eta)\rho)^2}{1-\rho ^{2}}}$, with $\frac{(1-(\alpha+\eta)\rho)^2}{1-\rho ^{2}}>1$ for $u$ close enough to $1$.

\noindent It remains to study $$\displaystyle\int_{1-(1-u)^{\alpha^2}}^{u}
\left(1-\Phi\left(\frac{\Phi^{-1}(u)-\rho\Phi^{-1}(v)}{\sqrt{1-\rho ^{2}}}\right)\right)dv.$$Recall that for $x>0$,  $\displaystyle1-\Phi(x)\leqslant\frac{e^{-\frac{x^2}{2}}}{\sqrt{2\pi}x}$. Thus we have 

\begin{align*}
&\int_{1-(1-u)^{\alpha^2}}^{u}
\left(1-\Phi\left(\frac{\Phi^{-1}(u)-\rho\Phi^{-1}(v)}{\sqrt{1-\rho ^{2}}}\right)\right)dv\\
&\leqslant\int_{1-(1-u)^{\alpha^2}}^{u} 
\frac{
e^{-\frac{1}{2}\frac{\Phi^{-1}(u)^2}{1-\rho ^{2}}}
e^{-\frac{1}{2}\frac{\rho^2\Phi^{-1}(v)^2}{1-\rho ^{2}}}
e^{\frac{\rho\Phi^{-1}(u)\Phi^{-1}(v)}{1-\rho ^{2}}}}
{\sqrt{2\pi}\frac{\Phi^{-1}(u)-\rho\Phi^{-1}(v)}{\sqrt{1-\rho ^{2}}}}dv\\
&\leqslant \frac{e^{-\frac{1}{2}\frac{\Phi^{-1}(u)^2}{1-\rho ^{2}}}}{\sqrt{2\pi}\Phi^{-1}(u)(1-\rho)}
\int_{1-(1-u)^{\alpha^2}}^{u} \sqrt{1-\rho ^{2}}
e^{-\frac{1}{2}\frac{\rho^2\Phi^{-1}(v)^2}{1-\rho ^{2}}}
e^{\frac{\rho\Phi^{-1}(u)\Phi^{-1}(v)}{1-\rho ^{2}}}dv\\
&=\frac{e^{-\frac{1}{2}\frac{\Phi^{-1}(u)^2}{1-\rho ^{2}}}}{\sqrt{2\pi}\Phi^{-1}(u)(1-\rho)}
\int_{\Phi^{-1}(1-(1-u)^{\alpha^2})}^{\Phi^{-1}(u)} \sqrt{1-\rho ^{2}}
e^{-\frac{1}{2}\frac{y^2}{1-\rho ^{2}}}
e^{\frac{\rho\Phi^{-1}(u)y}{1-\rho ^{2}}}\frac{dy}{\sqrt{2\pi}}\\
&=\frac{e^{-\frac{1}{2}\frac{\Phi^{-1}(u)^2}{1-\rho ^{2}}}}{\sqrt{2\pi}\Phi^{-1}(u)(1-\rho)}
\int_{\Phi^{-1}(1-(1-u)^{\alpha^2})-\rho\Phi^{-1}(u)}^{\Phi^{-1}(u)-\rho\Phi^{-1}(u)} \sqrt{1-\rho ^{2}}
e^{-\frac{1}{2}\frac{z^2}{1-\rho ^{2}}}
e^{\frac{1}{2}\frac{\rho^2\Phi^{-1}(u)^2}{1-\rho ^{2}}}\frac{dz}{\sqrt{2\pi}}\\
&=\frac{e^{-\frac{1}{2}\Phi^{-1}(u)^2}}{\sqrt{2\pi}\Phi^{-1}(u)(1-\rho)}
\int_{\Phi^{-1}(1-(1-u)^{\alpha^2})-\rho\Phi^{-1}(u)}^{\Phi^{-1}(u)-\rho\Phi^{-1}(u)} \sqrt{1-\rho ^{2}}
e^{-\frac{1}{2}\frac{z^2}{1-\rho ^{2}}}
\frac{dz}{\sqrt{2\pi}}=O\left(\frac{1-u}{\log(\frac{1}{1-u})}\right)
\end{align*}
since $\alpha<\frac{1}{\rho}(1-\sqrt{1-\rho ^{2}})<\rho$ and $\Phi^{-1}(1-(1-u)^{\alpha^2})\leqslant(\alpha+\eta)\Phi^{-1}(u)$, with $\eta$ arbitrarily small by choosing $u$ close to 1. Therefore this term is $O( (1-u)/\log(1/(1-u)) )$ near $1$.

\medskip
\noindent Now collecting the previous results, as $u\to 1$ we finally obtain
$$ \frac{u-C_\rho(u)}{h^2(u)}=O\left( \frac{1}{(1-u)\log^2(\frac{1}{1-u})}\right)$$which proves that 
it is integrable near $1$. By symmetry the same holds near $0$. We conclude that $(u-C_\rho(u))/h^2(u)$ is integrable on $(0,1)$.

\medskip

\medskip

\noindent{\bf Acknowledgements}{ We are grateful to Michel Ledoux who pointed out the question of the exact limiting constant in (\ref{BL}).}

\end{document}